%% file: uniqueness_paper_opttheoryapp_06.tex
\RequirePackage{fix-cm}
\documentclass[twocolumn]{svjour3}          
\smartqed  
\usepackage{graphicx}
\newcommand{\bb}{\mathbf}

\usepackage{verbatim} 

\usepackage{amssymb,amsmath}
\usepackage{tikz}
\usetikzlibrary{arrows}
\usepackage{pgfplots}

\usepackage{url} 

\begin{document}
	
	\title{Element-wise uniqueness, prior knowledge, and data-dependent resolution
	}
	
	
	\author{Keith Dillon         \and
		Yeshaiahu Fainman 
	}
	
	
	\institute{K. Dillon \at
		Tulane University, Department of Biomedical Engineering \\
		\email{kdillon1@tulane.edu}           
		\and
		Y. Fainman \at
        University of California, San Diego, Department of Electrical and Computer Engineering
	}
	
	\date{April 23, 2015}
	

	\maketitle
	
	\begin{abstract}
		Techniques for finding regularized solutions to underdetermined linear systems can be viewed as imposing prior knowledge on the unknown vector.
		The success of modern techniques, which can impose priors such as sparsity and non-negativity, is the result of advances in optimization algorithms to solve problems which lack closed-form solutions.
		Techniques for characterization and analysis of the system to determined when information is recoverable, however, still typically rely on closed-form solution techniques such as singular value decomposition or a filter cutoff, for example.
		In this letter we pose optimization approaches to broaden the approach to system characterization.
		We start by deriving conditions for when each unknown element of a system admits a unique solution, 
		subject to a broad class of types of prior knowledge.
		With this approach we can pose a convex optimization problem to find ``how unique'' each element of the solution is, which may be viewed as a generalization of resolution to incorporate prior knowledge.
		We find that the result varies with the unknown vector itself, i.e. is data-dependent, such as when the sparsity of the solution improves the chance it can be uniquely reconstructed. 
		The approach can be used to analyze systems on a case-by-case basis, estimate the amount of important information present in the data, and quantitatively understand the degree to which the regularized solution may be trusted. 
		
		\keywords{Underdetermined \and Regularization \and Resolution \and Super-resolution  \and Sparsity}
	\end{abstract}
	
	\section{Introduction}
	\label{intro}
	
	We focus on techniques that use norms such as the $\ell_1$-norm (sum of absolute elements) or the $\ell_\infty$-norm (maximum absolute element) for regularization and/or denoising of an underdetermined linear system, $\bb A \bb x = \bb b$, where $\bb A$ is a known $m \times n$ matrix with $m < n$, $\bb b$ is our known data vector, and $\bb x$ is the unknown vector we seek. 
	These techniques are generally not solvable in closed form (unlike for example, regularization with the $\ell_2$-norm).
	However modern optimization techniques can incorporate such information without difficulty using convex inequality constraints. 
	In this paper we will address approaches which may be specifically formulated as linear inequality constraints (see Appendix for some examples).
	In broad terms, instead of considering the restrictions on $\bb x$ in the set $\left\{ \bb x \;|\; \bb A \bb x = \bb b \right\}$,
	we focus on the (hopefully smaller and hence more informative) set  $\left\{ \bb x \;|\; \bb A \bb x = \bb b , \; \bb D \bb x \ge \bb d \right\}$.   
	We first ask the question: subject to this new constraint, has $\bb x$ become unique?
	We can extend this to the question: how well does this new information improve our ability to resolve $\bb x$?

	Uniqueness has been extensively studied for the case of $\ell_1$-regularization, where we are concerned for example with whether the solution found via Basis Pursuit \cite{chen_atomic_2001} is unique. 
	This is especially interesting because under the right conditions this solution is the optimally sparse solution (e.g. \cite{donoho_optimally_2003}). 
	Published conditions for $\ell_1$-regularized uniqueness come in several forms, such as the restricted isometry property \cite{candes_restricted_2008}, the null-space property \cite{donoho_optimally_2003}, and neighborliness properties \cite{donoho_neighborly_2005,donoho_neighborliness_2005}.
	Non-negativity constraints 
	have received increased interest recently due to their relationship to the above $\ell_1$-regularized case
	\cite{donoho_sparse_2005,bruckstein_uniqueness_2008}.
	In this application, if the true solution is sparse enough (and a necessary condition for the matrix holds), then the system has a unique non-negative solution.
	There is no regularization in this case, the non-negativity is directly applied as deterministic constraints on the solution.  
	Box constraints on $\bb x$ are a related case which has received some interest as well \cite{petra_3d_2009,donoho_counting_2010}.
	In all these approaches, however, the goal is a single cutoff which may be determined for the system itself. Whereas as we show in this letter, the answer actually varies, generally depending on the data and even between elements of the unknown vector.

	General views of resolution as relating to uniquely-determinable solutions have been proposed in inverse problem theory, in particular Backus-Gilbert theory \cite{bertero_linear_1985}, where the goal is to find a closed-form resolution estimate involving the $\ell_2$-norm. 
	Stark \cite{stark_generalizing_2008} proposed a broad framework to incorporate prior knowledge (and several other extensions), for which he proposed a general optimization approach, again to find a cutoff for the system resolution.
	A different direction is introduced by Cand\`es \cite{candes_towards_2013} for super-resolution, where gains due to sparsity of the unknown vector are described as yielding a super-resolution factor, essentially a higher resolution cutoff.

	In this paper we will formulate a novel approach to uniqueness by providing conditions on an element-wise (i.e. coordinate-wise) basis.
	The element-wise approach allows us to directly use convex optimization theory, as well as to make the relationship to the classical case (i.e. with no prior knowledge) clear.
	Further, we may relax the conditions with a test for uniqueness that, when it fails, can yield a kind of resolution estimate for the system. 
	Finally, we provide a simulation for a super-resolution scenario demonstrating how resolution varies intuitively with the choice of prior knowledge used and with the object itself.

\section{Methods} 

In our analysis we will neglect noise and model errors, presuming they are addressed by the denoising technique, if used, and so assume our underdetermined system $\bb A \bb x = \bb b$ has infinite solutions, which form the set, 
\begin{align}
F_{EC} =  \{ \bb x \in \mathbb{R}^n | \bb A \bb x = \bb b \}.
\label{ECset}
\end{align}
The subscript ``EC'' implies the solutions are 
equality-constrained.
In this paper we will consider the following set which has an added restriction representing our prior knowledge about the solution, 
\begin{align}
F_M =  \{ \bb x \in \mathbb{R}^n \; | \; \bb A \bb x = \bb b , \bb D \bb x \ge \bb d \},
\label{Fmixed}
\end{align}
where the subscript ``M'' implies mixed constraints.
By defining $\bb A$, $\bb D$, $\bb b$, and $\bb d$ in Eq. (\ref{Fmixed}) appropriately we may represent a variety of cases (see Appendix).
For example, we can consider the incorporation of non-negativity, as well as forms of regularization and denoising, and combinations of these.

\subsection{Uniqueness Conditions} 

Our first goal is to derive conditions for uniqueness of the $k^{th}$ element $x_k$ in $F_M$, for any selected $k \in \left\{1,...,n\right\}$.
To do this we will use optimization problems to solve for bounds on each element of $\bb x$.
When we refer to bounds on an element, we imply the maximum and minimum values that unknown element may take, given the information we have, as we investigated in \cite{dillon_bounding_2013}.
The bounds of the $k^{th}$ element (for any $k \in \left\{ 1, ..., n \right\}$) of a solution to a system are the scalar values given by
\begin{align}
  x_k^{(max)} &= \max \{ x_k \in \mathbb{R} | \bb x \in F \}, \\
  x_k^{(min)} &= \min \{ x_k \in \mathbb{R} | \bb x \in F \}.
\end{align}
Clearly an element $x_k$ is uniquely-determined if $x_k^{(max)} = x_k^{(min)}$.
We can test whether this is the case with the optimization problem,
\begin{align}
  \begin{array}{c}
    \delta_k = \\ 
   \; \\
   \; \\
   \; \\
   \;
  \end{array}
  \begin{array}{c}
    \underset{\bb x}{\max} \; x_k \\
    \bb A \bb x = \bb b \\
    \bb D \bb x \ge \bb d \\
    \; \\
    \;
  \end{array}
  \begin{array}{c}
    - \\ 
   \; \\
   \; \\
   \; \\
   \;
  \end{array}
  \begin{array}{c}
    \underset{\bb x}{\min} \; x_k \\
    \bb A \bb x = \bb b \\
    \bb D \bb x \ge \bb d \\
    \; \\
    \;
  \end{array}
    \begin{array}{c}
    = \\ 
   \; \\
   \; \\
   \; \\
   \;
  \end{array}
  \begin{array}{c}
    \underset{\bb x, \bb x'}{\max} \; ( x_k -x'_k )\\
    \bb A \bb x = \bb b \\
    \bb A \bb x' = \bb b \\
    \bb D \bb x \ge \bb d \\
    \bb D \bb x' \ge \bb d .    
  \end{array}
  \label{bounds_primal}
\end{align}
If the optimal value is $\delta_k=0$, then $x_k^{(max)} = x_k^{(min)}$ and $x_k$ must be uniquely determined.

Eqs. (\ref{bounds_primal}) form a linear program, and we can use duality theory for linear programming \cite{dantzig_linear_1998} to find an upper bound on $\delta_k$.
The dual of Eqs. (\ref{bounds_primal}) can be written as
\begin{align}
  \begin{array}{c}
    \tilde{\delta}_k = \\ 
   \; \\
   \; \\
   \; \\
   \;
  \end{array}
  \begin{array}{c}
    \underset{\bb y,\bb y', \bb z, \bb z'}{\min} \; 
    \bb b^T \left( \bb y - \bb y' \right) + \bb d^T \left( \bb z - \bb z' \right)\\
    \bb A^T \bb y + \bb D^T \bb z = \bb e_k \\
    \bb A^T \bb y' + \bb D^T \bb z' = -\bb e_k \\
    \bb z \le \bb 0 \\
    \bb z' \ge \bb 0. 
  \end{array}
\label{mixed_lp_gap_dual}
\end{align}
We form uniqueness conditions by requiring a feasible point exists such that the objective equals zero, giving the conditions,
\begin{align}
  \begin{array}{c}
    \bb b^T \left( \bb y - \bb y' \right) + \bb d^T \left( \bb z - \bb z' \right)  = 0\\
    \bb A^T \bb y + \bb D^T \bb z = \bb e_k \\
    \bb A^T \bb y' + \bb D^T \bb z' = -\bb e_k \\
    \bb z \le \bb 0 \\
    \bb z' \ge \bb 0. 
  \end{array}
\label{mixed_lp_gap_dual_cond}
\end{align}
As Eqs. (\ref{bounds_primal}) is the difference between a maximum and minimum over the same set, $\delta_k \ge 0$. 
Further, if a solution exists to Eqs. (\ref{mixed_lp_gap_dual_cond}) then $\tilde{\delta}_k = 0$, since $\delta_k \le \tilde{\delta}_k$, and by duality theory we have $0 \le \delta_k \le \tilde{\delta}_k=0$. 
Further, strong duality holds for linear programs under very general conditions (which we presume to hold), which means $\delta_k = \tilde{\delta}_k$.

To understand the conditions of Eqs. (\ref{mixed_lp_gap_dual_cond}), note that if $\bb D$ and $\bb d$ are set to zeros (and hence we are back to the classical case of Eq. (ref{ECset})), then the conditions can be met for any $\bb y$ such that $\bb A^T \bb y = \bb e_k$. 
So for the classical case $\bb y$ is simply a (transposed) row of the left inverse of $\bb A$,  
and this condition can be viewed as an element-wise version of the condition that $\bb A$ is non-singular. 
Note that this classical condition does not depend on $\bb b$, while the conditions of Eqs. (\ref{mixed_lp_gap_dual_cond}) do.
Since $\bb b = \bb A \bb x$, uniqueness when there is prior knowledge included will (in general) depend on the particular value of $\bb x$ in each case. 
Further, for the case where there is no solution to $\bb A^T \bb y = \bb e_k$, we may still able able to solve the equation $\bb A^T \bb y + \bb D^T \bb z = \bb e_k$, for example, with an appropriate choice of $\bb z$.
So the prior knowledge represented by $\bb D \bb x \ge \bb d$ results in a restriction on the possible $\bb x$, and a  relaxation of the uniqueness conditions. 

\subsection{Resolution } 

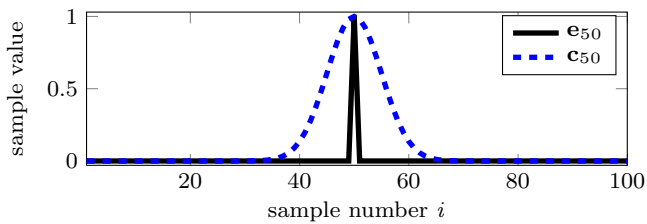
\begin{figure}[!t] \centering 
    \input{SITL07_genericreso_101.tex} 
    \caption{$\bb e_k$ versus relaxed result for $k=50$ with  $n=100$.}
    \label{101genericreso}
\end{figure}

Now we will consider a way to relax the uniqueness conditions which provides a metric for the system, which we can then use to compare the improvement due to various cases of prior knowledge.
To motivate the approach, consider the non-bounded case again.
If a ${\bb y}$ can be found such that $\bb A^T {\bb y} = \bb e_k$, then we can compute ${\bb y}^T \bb b =  {\bb y}^T \bb A \bb x = \bb e_k^T \bb x = x_k$. 
So ${\bb y}$ is a linear functional that computes $x_k$.
In the event that finding such a functional is not possible, our goal is to find one that gets as close as possible.
As depicted in Fig. \ref{101genericreso}, we replace $\bb e_k$ with a vector $\bb c$ that has some spread over multiple elements.
To find the $\bb c$ closest to $\bb e_k$ we use an optimization problem such as the following,
\begin{align}
  \begin{array}{c}
    d^{(EC)}_k = \\ 
   \; \\
   \; \\
   \;
  \end{array}
  \begin{array}{c}
     \underset{\bb c, \bb y}{\min} \; \Vert \bb c \Vert \\
    \bb A^T \bb y = \bb c \\
    \bb c \ge 0 \\
    c_k = 1.
  \end{array}
\label{reso_equality-constrained}
\end{align}

In the case where $\bb A^T \bb y = \bb e_k$ has a solution, Eqs. (\ref{reso_equality-constrained}) will achieve $\bb c = \bb e_k$.
Otherwise, the result will be a metric of how similar $\bb c$ could be made to $\bb e_k$, depending on the choice of norm in the objective. 
To provide a intuitively-meaningful metric, we included the constraint $\bb c \ge \bb 0$ and for the norm use a $\ell_2$-norm weighted with distance (in terms of spatial or temporal location of the samples) from the $k^{th}$ element.
If the weighting increases quadratically, $\bb c$ can be viewed as a distribution and the metric can be viewed as its variance.
So the optimization seeks the distribution about the element of interest with the minimum variance, which may be uniquely-determined.

Similarly, the conditions of Eqs. (\ref{mixed_lp_gap_dual_cond}) give the optimization problem subject to prior knowledge,
\begin{align}
\begin{array}{c}
d^{(M)}_k = \\ 
\; \\
\; \\
\; \\
\; \\
\; \\
\; \\
\;
\end{array}
\begin{array}{c}
 \underset{\bb c, \bb y, \bb y', \bb z, \bb z'}{\min} \; \Vert \bb c \Vert \\
    \bb b^T \left( \bb y - \bb y' \right) + \bb d^T \left( \bb z - \bb z' \right)  = 0\\
    \bb A^T \bb y + \bb D^T \bb z = \bb e_k \\
    \bb A^T \bb y' + \bb D^T \bb z' = -\bb e_k \\
    \bb z \le \bb 0 \\
    \bb z' \ge \bb 0 \\
\bb c \ge 0 \\
c_k = 1.
\end{array}
\label{reso_mixed}
\end{align}
The constraints are linear so this is a convex optimization problem.

\section{Simulation} 

\begin{figure}[!t] \centering 
    \input{go1DSim15_kernel_xtrue_b.tex} 
    \caption{Test input $\bb x^{(true)}$, the true values of the unknown vector, kernel convolved with $\bb x$ prior to downsampling, and measured data $\bb b = \bb A \bb x^{(true)}$.}
    \label{101kernel}
\end{figure}
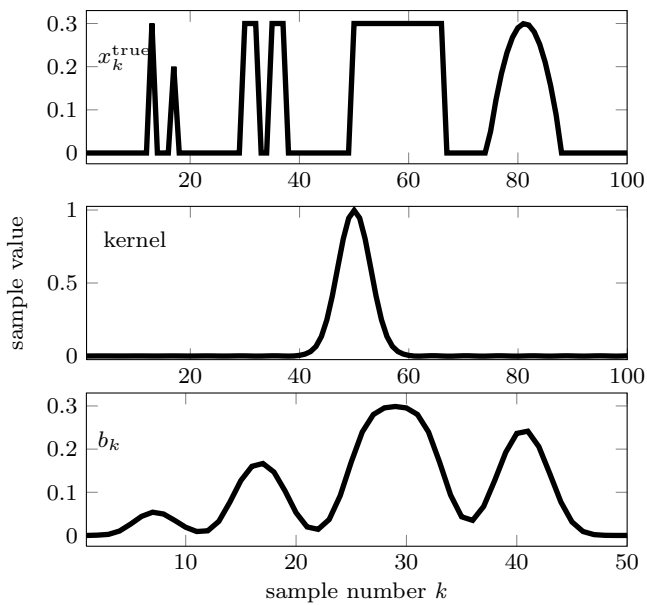

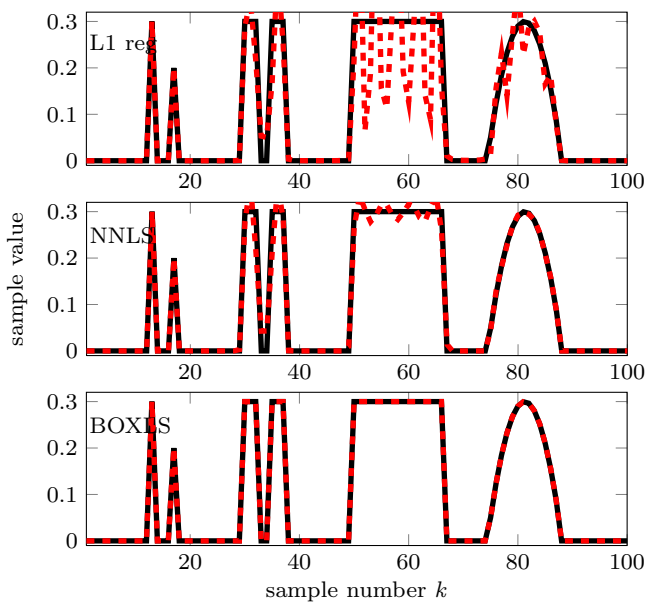
\begin{figure}[!t] \centering 
    \input{go1DSim16_stateart_noiseless.tex} 
    \caption{Regularized estimates with different techniques: L1-regularized, a.k.a Basis Pursuit (L1); Non-negative Least-Squares (NNLS);  Box-constrained Least-Squares (BOXLS). Dashed red trace is estimate, solid black trace is true $\bb x$ for comparison.}
    \label{101stateart}
\end{figure}

To demonstrate the resolution estimate we simulated a one-dimensional system which lowpass filters and downsamples the input.
The matrix $\bb A$ convolves the input vector with a Gaussian-shaped kernel then downsamples the result by a factor of two. So $\bb A$ is $m \times n$ with $n=100$ and $m=50$.
The true $\bb x$, convolution kernel, and lowpass version are shown in Fig. \ref{101kernel}. 
In Fig. \ref{101stateart} we show $\bb x^{(true)}$ versus some common regularized estimates.
In addition to the conventional Basis Pursuit and non-negative techniques, we also used a technique analogous to non-negative least-squares but with box constraints (both a lower and upper constraint) on $\bb x$.
Specifically, we computed the minimum 2-norm subject to the constraint $0 \le x_i \le 0.3$ for each element of $\bb x$. 

We see that $\ell_1$-regularization did not yield a very good result; on the left side of $\bb x$, where the signal is locally sparse (the number of nonzero values in a local intervals is low) we have a correct estimate, but on the right side of the plot where $\bb x$ is denser, the estimate becomes incorrect.
The non-negative least-squares estimate gave a better result, but still was incorrect in the densest region in the center-right of the plot.
The box-constrained least-squares produced an apparently-perfect result, as we used both the true upper and lower limits as prior knowledge.

\begin{figure}[!t] \centering 
    \input{go1DSim16_reso.tex} 
    \caption{Resolution estimate for each sample for different cases.}
    \label{101reso}
\end{figure}
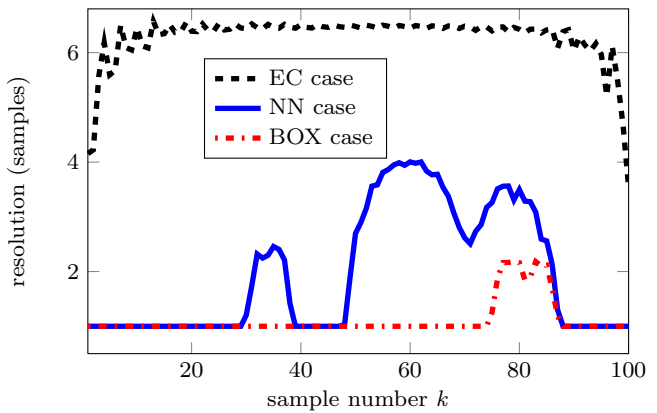

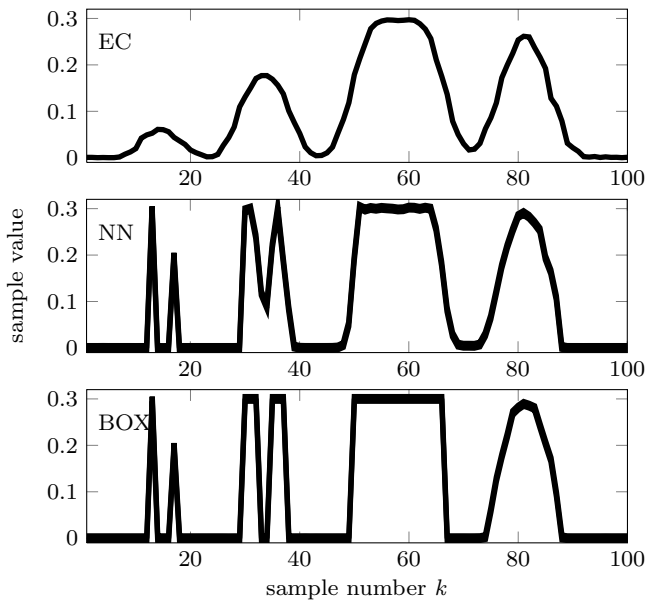
\begin{figure}[!t] \centering 
    \input{go1DSim16_lowres.tex} 
    \caption{Low resolution estimates of $\bb x$ for different cases.}
    \label{101lowres}
\end{figure}

We used CVX \cite{grant_cvx:_2012,grant_graph_2008} to solve the optimization problems.
%
The resolution functions $\bb c$ were estimated for each element of $\bb x$ for three different cases, equality-constrained, non-negative, and box-constrained (see Appendix for formulation as $\bb D$ and $\bb d$).
The element-wise resolution estimates for the different cases are plotted in Fig. \ref{101reso},
and ``Low-resolution'' estimates using the resolution cells, i.e. $\bb c^T \bb b$, analogous to $\bb e_k^T \bb x$, are plotted in Fig. \ref{101lowres}. 
We see that the equality-constrained case returns an essentially constant resolution (except for edge effects) which
approximates the amount of low-pass filtering performed by the kernel.
The box-constrained case achieves best resolution (one sample, implying $\bb c = \bb e_k$) for most of the elements, as we might have guessed given the accurate reconstruction, except in a small interval around sample 80.
This poorer-resolution region underlines the fact that an apparently-perfect regularized reconstruction does not necessarily imply a unique solution and hence a sufficient system resolution. 
The non-negative case achieved results in between. 
It did achieve maximum resolution for regions where the signal was locally the sparsest, i.e. the pair of spikes and the regions of zeros.

The data-dependent performance for the constrained cases can be intuitively understood in terms of active constraints (e.g., an inequality constraint such as $x_k \ge 0$ holds with equality, $x_k=0$).
Recall our data, $\bb b$, is a collection of mixtures of elements of $\bb x$ about local regions, as defined by the shape of the kernel. 
When several of those elements are equal to the constraint, the range of possible combinations gets restricted and the resolution potentially improves.
For non-negativity constraints, active constraints occur at elements taking the value zero, hence we see the effect of sparsity on resolution.
Similarly for box constraints we find the resulting resolution is best for regions where the true signal is either zero or 0.3.

	\section{Discussion} 
	\raggedbottom 
	In this paper we gave conditions for each element of an underdetermined linear system of equations and inequality constraints to be uniquely determined.
	We used element-wise uniqueness as the basis for a definition of resolution, as the size of the smallest resolution cell that can be found.
	This allowed us to compute resolution on an element-wise basis subject to prior knowledge in the form of inequality constraints.

	The resolution cell estimate is very interesting when inequality constraints are included, as it yields a data-dependent result. 
	In the case of non-negativity, this result depends on the sparsity of the elements which are mixed with our element of interest.
	In the case of more general inequality constraints, the sparsity condition would be replaced with a metric of the number of active constraints. 
	For the simulated case, essentially a super-resolution problem, this mixing is localized so we see the effect due to the active constraints in local regions.
	For such a system, concentrated resolution estimate makes sense. 
	For more arbitrary systems a concentrated resolution cell may not be achievable.
	This would imply that the ambiguity between high-resolution elements cannot be explained with any locally-concentrated combination. 
	
	There are a variety of ways one could estimate the most compact resolution cell for each pixel.
	The $\ell_2$-norm was used here as it yielded an intuitive interpretation in terms of the variance of a distribution over space or time.

	The computational complexity of the technique requires one optimization problem per estimate.
	When the low-resolution distributions are large, the estimates at neighboring elements, are mostly redundant. 
	We can therefore choose to increase the spacing between estimates such that we still achieve a covering of all elements. 
	This may follow a variety of strategies, for example we may randomly choose uncovered elements in a greedy approach or simply step across the elements with a stepsize equal to the size of the previous resolution cell.

	\section{Appendix}
	
	In this appendix we will describe how a selection of variations on prior knowledge can be formulated as linear inequality constraints.
	First, the classical case, with no prior knowledge is based on the solution set,
	\begin{align}
	F_{EC} =  \{ \bb x \in \mathbb{R}^n | \bb A \bb x = \bb b \}.
	\label{ECset}
	\end{align}
	Application of our bounds testing problem with the feasible set $\bb x \in F_{EC}$ forms an equality-constrained linear program \cite{gill_numerical_1991}, for which the optimality conditions give the rowspace condition $\bb A^T \bb y = \bb e_k$.
	
	%
	%

	\subsection{Non-negativity} 
	The simplest restriction on $\bb x$ is non-negativity, resulting in the solution set
	\begin{align}
	F_{NN} =  \{ \bb x \in \mathbb{R}^n | \bb A \bb x = \bb b , \bb x \ge \bb 0 \}.
	\end{align}
	This can be implemented in our system with the simple definitions, $\bb D = \bb I$, $\bb d = \bb 0$, using the identity matrix and a vector of zeros.

	\subsection{Basis Pursuit} 
	
	Basis pursuit \cite{chen_atomic_2001} is a regularization technique that utilizes the $\ell_1$-norm.
	We are interested in whether we have a unique optimal solution to the Basis Pursuit problem, 
	\begin{align}
	\begin{array}{c}
	\alpha = \\ 
	\;
	\end{array}
	\begin{array}{c}
	\underset{\bb x}{\min} \; \Vert \bb x \Vert_1 \\
	\bb A \bb x = \bb b .
	\end{array}
	\label{basis_pursuit}
	\end{align}
	This can be tested by analyzing in the uniqueness of the solutions in the following set,
	\begin{align}
	F_{BP} &=  \{ \bb x \in \mathbb{R}^n | \bb A \bb x = \bb b , \, \Vert \bb x \Vert_1 = \alpha \} \notag \\
	&=  \{ \bb x \in \mathbb{R}^n | \bb A \bb x = \bb b , \, \Vert \bb x \Vert_1 \le \alpha \}.
	\label{Fbp}
	\end{align}
	This is equivalent to the following non-negative system,
	\begin{align}
	F_{NN} =  \{  \hat{\bb x} \in \mathbb{R}^{2n} | \hat{\bb A}  \hat{\bb x} = \hat{\bb b} ,  \hat{\bb x} \ge \bb 0 \}.
	\label{Fbpnn}
	\end{align}
	With the definitions
	\begin{align}
	\hat{\bb A} = 
	\begin{pmatrix} 
	\bb A, -\bb A \\ \bb 1^T 
	\end{pmatrix}
	\end{align}
	\begin{align}
	\hat{\bb b} = 
	\begin{pmatrix} 
	\bb b \\ \alpha.
	\end{pmatrix}
	\end{align}
	This can be seen by defining $\bb x = \hat{\bb x}_{(1)} - \hat{\bb x}_{(2)}$, 
	where 
	$\hat{\bb x}^T = \begin{pmatrix}\hat{\bb x}_{(1)}^T, \hat{\bb x}_{(2)}^T\end{pmatrix}$ 
	and $\hat{\bb x}_{(1)} \ge \bb 0$, $\hat{\bb x}_{(2)}\ge \bb 0$. 
	We relate bounds found using the feasible set of Eq. (\ref{Fbpnn}) to the bounds for the set of Eq. (\ref{Fbp}) by noting that at the minimum, where we get $\alpha$ as the optimal for Eq. (\ref{basis_pursuit}), $\hat{\bb x}_{(1)}$ and $\hat{\bb x}_{(1)}$ are complementary. 
	If they were not, we could take advantage of this fact to reduce the minimum of $\Vert \bb x \Vert_1 = \hat{\bb x}_{(1)} + \hat{\bb x}_{(2)}$ further.

	\subsection{Box Constrained Regularization} 

	A similar example is box constraints, defined as the system
	\begin{align}
	F_{BOX} =  \{ \bb x \in \mathbb{R}^n \; | \; \bb A \bb x = \bb b , \; \bb d_{min} \le \bb x \le \bb d_{max} \}.
	\label{Fbox}
	\end{align}
	Here $\bb d_{min}$ and  $\bb d_{max}$ are vectors defining the box.
	We can formulate this as Eq. (\ref{Fmixed}) with the definitions
	\begin{align}
	\bb D = 
	\begin{pmatrix} 
	+\bb I \\ -\bb I .
	\end{pmatrix},
	\end{align}
	\begin{align}
	\bb d = 
	\begin{pmatrix} 
	+\bb d_{min} \\ -\bb d_{max}.
	\end{pmatrix}.
	\end{align}
	We can view this as a case of regularization with the infinity norm, e.g., 
	\begin{align}
	F_{BOX} =  \{ \bb x \in \mathbb{R}^n | \bb A \bb x = \bb b , \, \Vert \bb x \Vert_\infty \le d \}.
	\label{Fbp}
	\end{align}

	\subsection{Denoising} 

	Denoising can be viewed as a dual to regularization, where rather than require the solution set be regular, we require the error in the linear system to be regular, as in the following, 
	\begin{align}
	F_{DN} =  \{ \bb x \in \mathbb{R}^n \; | \; \Vert \bb A \bb x - \bb b \Vert \le \sigma,  \}.
	\label{Fsigma}
	\end{align}
	Where $\sigma$ is the minimum over $\Vert \bb A \bb x - \bb b \Vert$. 
	We can form a denoised version of the non-negativity case using the infinity norm as follows,
	\begin{align}
	F_{NND} =  \{ \bb x \in \mathbb{R}^n \; | \; \Vert \bb A \bb x - \bb b \Vert_\infty \le \sigma_{NN}, \; \bb x \ge \bb 0 \}.
	\label{Fsigma}
	\end{align}
	This can be formulated as mixed constraints with no linear constraint term 
	(i.e. ``$\bb A$'' and ``$\bb b$'' in the original linear system are all zeros), 
	and with 
	\begin{align}
	\bb D = 
	\begin{pmatrix}
	-\bb A \\
	\bb A \\
	\bb I
	\end{pmatrix},
	\end{align}
	\begin{align}
	\bb d = 
	\begin{pmatrix}
	-\bb b - \sigma_{NN} \bb 1 \\
	\bb b - \sigma_{NN} \bb 1 \\
	\bb 0
	\end{pmatrix}.
	\end{align}

	
	\bibliographystyle{spmpsci_modified}      

	\bibliography{zoterorefs}   
	
	
\end{document}

%% file: SITL07_genericreso_101.tex
%
%
\begin{tikzpicture}[font=\small]

\begin{axis}[%
width=2.8in,
height=0.8in,
ylabel style={yshift=-1em},
xlabel style={yshift=.5em},
scale only axis,
xmin=1,
xmax=100,
xlabel={$\text{sample number }{\it i}$},
ymin=-0.03,
ymax=1.03,
ylabel={sample value},
legend style={draw=black,fill=white,legend cell align=left}
]
\addplot [
color=black,
solid,
line width=2.0pt
]
table[row sep=crcr]{
1 0\\
2 0\\
3 0\\
4 0\\
5 0\\
6 0\\
7 0\\
8 0\\
9 0\\
10 0\\
11 0\\
12 0\\
13 0\\
14 0\\
15 0\\
16 0\\
17 0\\
18 0\\
19 0\\
20 0\\
21 0\\
22 0\\
23 0\\
24 0\\
25 0\\
26 0\\
27 0\\
28 0\\
29 0\\
30 0\\
31 0\\
32 0\\
33 0\\
34 0\\
35 0\\
36 0\\
37 0\\
38 0\\
39 0\\
40 0\\
41 0\\
42 0\\
43 0\\
44 0\\
45 0\\
46 0\\
47 0\\
48 0\\
49 0\\
50 1\\
51 0\\
52 0\\
53 0\\
54 0\\
55 0\\
56 0\\
57 0\\
58 0\\
59 0\\
60 0\\
61 0\\
62 0\\
63 0\\
64 0\\
65 0\\
66 0\\
67 0\\
68 0\\
69 0\\
70 0\\
71 0\\
72 0\\
73 0\\
74 0\\
75 0\\
76 0\\
77 0\\
78 0\\
79 0\\
80 0\\
81 0\\
82 0\\
83 0\\
84 0\\
85 0\\
86 0\\
87 0\\
88 0\\
89 0\\
90 0\\
91 0\\
92 0\\
93 0\\
94 0\\
95 0\\
96 0\\
97 0\\
98 0\\
99 0\\
100 0\\
};
\addlegendentry{${\bf e}_{\text{50}}$};

\addplot [
color=blue,
dashed,
line width=2.0pt
]
table[row sep=crcr]{
1 1.3969439431471e-21\\
2 9.72098502030078e-21\\
3 6.49934797207089e-20\\
4 4.17501005585051e-19\\
5 2.57675710915498e-18\\
6 1.5279799682873e-17\\
7 8.70542662229619e-17\\
8 4.76530473529909e-16\\
9 2.50622188714527e-15\\
10 1.26641655490942e-14\\
11 6.14839641270476e-14\\
12 2.8679750088881e-13\\
13 1.28533722513365e-12\\
14 5.53461007170101e-12\\
15 2.28973484564555e-11\\
16 9.10147076448794e-11\\
17 3.47589128123992e-10\\
18 1.27540762952604e-09\\
19 4.49634946228087e-09\\
20 1.52299797447126e-08\\
21 4.9564053191725e-08\\
22 1.5497531357029e-07\\
23 4.65571571578309e-07\\
24 1.34381227763152e-06\\
25 3.72665317207867e-06\\
26 9.92950430585108e-06\\
27 2.54193465161992e-05\\
28 6.25215037748203e-05\\
29 0.000147748360232034\\
30 0.000335462627902512\\
31 0.000731802418880473\\
32 0.00153381067932446\\
33 0.00308871540823677\\
34 0.00597602289500594\\
35 0.0111089965382423\\
36 0.0198410947443703\\
37 0.0340474547345993\\
38 0.0561347628341337\\
39 0.0889216174593863\\
40 0.135335283236613\\
41 0.197898699083615\\
42 0.278037300453194\\
43 0.3753110988514\\
44 0.486752255959972\\
45 0.606530659712633\\
46 0.726149037073691\\
47 0.835270211411272\\
48 0.923116346386636\\
49 0.980198673306755\\
50 1\\
51 0.980198673306755\\
52 0.923116346386636\\
53 0.835270211411272\\
54 0.726149037073691\\
55 0.606530659712633\\
56 0.486752255959972\\
57 0.3753110988514\\
58 0.278037300453194\\
59 0.197898699083615\\
60 0.135335283236613\\
61 0.0889216174593863\\
62 0.0561347628341337\\
63 0.0340474547345993\\
64 0.0198410947443703\\
65 0.0111089965382423\\
66 0.00597602289500594\\
67 0.00308871540823677\\
68 0.00153381067932446\\
69 0.000731802418880473\\
70 0.000335462627902512\\
71 0.000147748360232034\\
72 6.25215037748203e-05\\
73 2.54193465161992e-05\\
74 9.92950430585108e-06\\
75 3.72665317207867e-06\\
76 1.34381227763152e-06\\
77 4.65571571578309e-07\\
78 1.5497531357029e-07\\
79 4.9564053191725e-08\\
80 1.52299797447126e-08\\
81 4.49634946228087e-09\\
82 1.27540762952604e-09\\
83 3.47589128123992e-10\\
84 9.10147076448794e-11\\
85 2.28973484564555e-11\\
86 5.53461007170101e-12\\
87 1.28533722513365e-12\\
88 2.8679750088881e-13\\
89 6.14839641270476e-14\\
90 1.26641655490942e-14\\
91 2.50622188714527e-15\\
92 4.76530473529909e-16\\
93 8.70542662229619e-17\\
94 1.5279799682873e-17\\
95 2.57675710915498e-18\\
96 4.17501005585051e-19\\
97 6.49934797207089e-20\\
98 9.72098502030078e-21\\
99 1.3969439431471e-21\\
100 1.92874984796392e-22\\
};
\addlegendentry{${\bf c}_{\text{50}}$};

\end{axis}
\end{tikzpicture}%

%% file: go1DSim15_kernel_xtrue_b.tex
%
%
\begin{tikzpicture}[font=\small]

\begin{axis}[%
width=2.8in,
height=.8in,
ylabel style={yshift=-1em},
xlabel style={yshift=.5em},
scale only axis,
xmin=1,
xmax=100,
ymin=-0.025,
ymax=1.025,
ylabel={sample value},
name=plot2
]
\addplot [
color=black,
solid,
line width=2.0pt,
forget plot
]
table[row sep=crcr]{
1 -0.000201170010674827\\
2 -3.63836385232791e-05\\
3 0.000119888217457321\\
4 0.000175441281391935\\
5 8.87955443113691e-05\\
6 -8.11110744464101e-05\\
7 -0.000198868314466129\\
8 -0.000153462862863733\\
9 3.71991724241455e-05\\
10 0.000218460507946919\\
11 0.000219507089559718\\
12 1.34016048557295e-05\\
13 -0.000232609940358905\\
14 -0.000287541237298168\\
15 -7.18572080838552e-05\\
16 0.000241378505119733\\
17 0.000359012121676126\\
18 0.000139589707693277\\
19 -0.000244920515158371\\
20 -0.000435712024279488\\
21 -0.000218291385399901\\
22 0.00024360921916899\\
23 0.000520093996458041\\
24 0.000310057156621264\\
25 -0.00023836035679411\\
26 -0.000615800691113354\\
27 -0.000417495852236953\\
28 0.000231245589246113\\
29 0.000728446878285163\\
30 0.000543629451599134\\
31 -0.000226754023900614\\
32 -0.00086671061441509\\
33 -0.000690948239471856\\
34 0.000235006344113622\\
35 0.00104727237440417\\
36 0.000877644345640801\\
37 -0.000189301828120397\\
38 -0.000939821947264449\\
39 0.000168608271501981\\
40 0.00424227184226158\\
41 0.0126820869517868\\
42 0.0297469479430173\\
43 0.0651388892453652\\
44 0.133417999573045\\
45 0.248129372616454\\
46 0.412048487141562\\
47 0.608756129960318\\
48 0.801818099211168\\
49 0.944600639479259\\
50 0.997627519428863\\
51 0.945202203816483\\
52 0.80245288489335\\
53 0.608824012810122\\
54 0.411484503410794\\
55 0.247465839830361\\
56 0.133282462878486\\
57 0.0656610482732304\\
58 0.0304346419862821\\
59 0.0128848203447144\\
60 0.00376606445833615\\
61 -0.000538550360104377\\
62 -0.00120906924486654\\
63 0.000236942859210039\\
64 0.00159946399302462\\
65 0.00138212797654622\\
66 -0.000137380105842968\\
67 -0.00142251944374893\\
68 -0.00126605170232502\\
69 8.79907496782349e-05\\
70 0.00127994166666154\\
71 0.00119093952997134\\
72 -2.21832860198456e-05\\
73 -0.00115344116836612\\
74 -0.00113991114197352\\
75 -4.98160765281035e-05\\
76 0.00104044303571012\\
77 0.00110410324753748\\
78 0.000123411289733483\\
79 -0.000937865209205865\\
80 -0.00107774300295394\\
81 -0.000196408415425361\\
82 0.000843101998676394\\
83 0.00105707380577112\\
84 0.000267717070547763\\
85 -0.000754205253291019\\
86 -0.00103962505886508\\
87 -0.000336658415786817\\
88 0.000669988044956866\\
89 0.00102387076644763\\
90 0.00040249660862301\\
91 -0.000590643290516943\\
92 -0.00100989251671374\\
93 -0.000463766575921327\\
94 0.000521842486653821\\
95 0.00100505169070802\\
96 0.000513781075318019\\
97 -0.000505263416953922\\
98 -0.0010639455435007\\
99 -0.00046228105756667\\
100 0.00116885455624571\\
};
\node[right, inner sep=0mm, text=black]
at (axis cs:4,0.8,0) {kernel};
\end{axis}

\begin{axis}[%
width=2.8in,
height=.8in,
ylabel style={yshift=-1em},
xlabel style={yshift=.5em},
scale only axis,
xmin=1,
xmax=100,
ymin=-0.025,
ymax=0.33,
at=(plot2.above north west),
anchor=below south west
]
\addplot [
color=black,
solid,
line width=2.0pt,
forget plot
]
table[row sep=crcr]{
1 0\\
2 0\\
3 0\\
4 0\\
5 0\\
6 0\\
7 0\\
8 0\\
9 0\\
10 0\\
11 0\\
12 0\\
13 0.3\\
14 0\\
15 0\\
16 0\\
17 0.2\\
18 0\\
19 0\\
20 0\\
21 0\\
22 0\\
23 0\\
24 0\\
25 0\\
26 0\\
27 0\\
28 0\\
29 0\\
30 0.3\\
31 0.3\\
32 0.3\\
33 0\\
34 0\\
35 0.3\\
36 0.3\\
37 0.3\\
38 0\\
39 0\\
40 0\\
41 0\\
42 0\\
43 0\\
44 0\\
45 0\\
46 0\\
47 0\\
48 0\\
49 0\\
50 0.3\\
51 0.3\\
52 0.3\\
53 0.3\\
54 0.3\\
55 0.3\\
56 0.3\\
57 0.3\\
58 0.3\\
59 0.3\\
60 0.3\\
61 0.3\\
62 0.3\\
63 0.3\\
64 0.3\\
65 0.3\\
66 0.3\\
67 0\\
68 0\\
69 0\\
70 0\\
71 0\\
72 0\\
73 0\\
74 0\\
75 0.05\\
76 0.1236\\
77 0.1844\\
78 0.2324\\
79 0.2676\\
80 0.29\\
81 0.2996\\
82 0.2964\\
83 0.2804\\
84 0.2516\\
85 0.21\\
86 0.1556\\
87 0.0884\\
88 0\\
89 0\\
90 0\\
91 0\\
92 0\\
93 0\\
94 0\\
95 0\\
96 0\\
97 0\\
98 0\\
99 0\\
100 0\\
};
\node[right, inner sep=0mm, text=black]
at (axis cs:3,0.22,0) {$ x^{\text{true}}_k$};
\end{axis}

\begin{axis}[%
width=2.8in,
height=.8in,
ylabel style={yshift=-1em},
xlabel style={yshift=.5em},
scale only axis,
xmin=1,
xmax=50,
xlabel={$\text{sample number }{\it{ k}}$},
ymin=-0.025,
ymax=0.33,
at=(plot2.below south west),
anchor=above north west
]
\addplot [
color=black,
solid,
line width=2.0pt,
forget plot
]
table[row sep=crcr]{
1 -6.12694964071085e-05\\
2 0.000681091850664482\\
3 0.0024593275083735\\
4 0.0102124694667635\\
5 0.0261837664147427\\
6 0.0441530503326892\\
7 0.0538422063902901\\
8 0.0496081981801179\\
9 0.0348962013909568\\
10 0.0189382330665687\\
11 0.0090833989130169\\
12 0.0107832673475488\\
13 0.0326509533402593\\
14 0.077084016802149\\
15 0.127250650900335\\
16 0.160135740803329\\
17 0.166704195761058\\
18 0.146787804178023\\
19 0.103156147045847\\
20 0.0526873738515859\\
21 0.0197172695559802\\
22 0.0141553574497943\\
23 0.0369330373924356\\
24 0.0921924060343505\\
25 0.170011540264942\\
26 0.239733735312132\\
27 0.280073672374224\\
28 0.295537370738014\\
29 0.298779649598472\\
30 0.295369962921334\\
31 0.280254999691765\\
32 0.239728656242474\\
33 0.170017599493876\\
34 0.093704266732369\\
35 0.0432100351278539\\
36 0.0349426075164734\\
37 0.06726235060794\\
38 0.127442140780142\\
39 0.192579510220374\\
40 0.236398150026008\\
41 0.241522056719297\\
42 0.205628959071169\\
43 0.142522578452307\\
44 0.0773215679161359\\
45 0.0313534825491226\\
46 0.00911625843066371\\
47 0.0018954672835177\\
48 0.000270127209610635\\
49 -1.26932858721166e-05\\
50 2.56127781357738e-05\\
};
\node[right, inner sep=0mm, text=black]
at (axis cs:2,0.22,0) {$b_k$};
\end{axis}
\end{tikzpicture}%

%% file: go1DSim16_stateart_noiseless.tex
%
%
\begin{tikzpicture}[font=\small]

\begin{axis}[%
width=2.8in,
height=.8in,
ylabel style={yshift=-1em},
xlabel style={yshift=.5em},
scale only axis,
xmin=1,
xmax=100,
ymin=-0.01,
ymax=0.32,
ylabel={sample value},
name=plot2
]
\addplot [
color=black,
solid,
line width=2.0pt,
forget plot
]
table[row sep=crcr]{
1 0\\
2 0\\
3 0\\
4 0\\
5 0\\
6 0\\
7 0\\
8 0\\
9 0\\
10 0\\
11 0\\
12 0\\
13 0.3\\
14 0\\
15 0\\
16 0\\
17 0.2\\
18 0\\
19 0\\
20 0\\
21 0\\
22 0\\
23 0\\
24 0\\
25 0\\
26 0\\
27 0\\
28 0\\
29 0\\
30 0.3\\
31 0.3\\
32 0.3\\
33 0\\
34 0\\
35 0.3\\
36 0.3\\
37 0.3\\
38 0\\
39 0\\
40 0\\
41 0\\
42 0\\
43 0\\
44 0\\
45 0\\
46 0\\
47 0\\
48 0\\
49 0\\
50 0.3\\
51 0.3\\
52 0.3\\
53 0.3\\
54 0.3\\
55 0.3\\
56 0.3\\
57 0.3\\
58 0.3\\
59 0.3\\
60 0.3\\
61 0.3\\
62 0.3\\
63 0.3\\
64 0.3\\
65 0.3\\
66 0.3\\
67 0\\
68 0\\
69 0\\
70 0\\
71 0\\
72 0\\
73 0\\
74 0\\
75 0.05\\
76 0.1236\\
77 0.1844\\
78 0.2324\\
79 0.2676\\
80 0.29\\
81 0.2996\\
82 0.2964\\
83 0.2804\\
84 0.2516\\
85 0.21\\
86 0.1556\\
87 0.0884\\
88 0\\
89 0\\
90 0\\
91 0\\
92 0\\
93 0\\
94 0\\
95 0\\
96 0\\
97 0\\
98 0\\
99 0\\
100 0\\
};
\addplot [
color=red,
dashed,
line width=2.0pt,
forget plot
]
table[row sep=crcr]{
1 -1.41383505116238e-09\\
2 -7.28589177852739e-09\\
3 -1.03046066809188e-08\\
4 -1.17539274880844e-08\\
5 -1.21152679092055e-08\\
6 -1.14502041264795e-08\\
7 -9.42647324152447e-09\\
8 -4.99409560949925e-09\\
9 4.81753943286866e-09\\
10 2.95916092797074e-08\\
11 1.10661942274888e-07\\
12 5.99365511184175e-07\\
13 0.299981048323361\\
14 2.87510053934092e-05\\
15 3.18581820777307e-06\\
16 4.62520350864574e-06\\
17 0.199951852702711\\
18 2.79400801402514e-05\\
19 1.38049760860359e-06\\
20 5.98751703170367e-07\\
21 4.21995782733402e-07\\
22 4.07979299574335e-07\\
23 5.17375114636542e-07\\
24 8.51357602916983e-07\\
25 1.77564576645198e-06\\
26 4.02922047236966e-06\\
27 7.21416748222512e-06\\
28 1.2676770553059e-05\\
29 0.00289818561402204\\
30 0.278982701116574\\
31 0.360175667890196\\
32 0.215593932240227\\
33 0.0425638416677359\\
34 0.0430149736272668\\
35 0.21420709639805\\
36 0.357803437583599\\
37 0.284227782222288\\
38 0.000400515462383474\\
39 1.23496314349403e-05\\
40 7.9932909313228e-06\\
41 8.69860505505331e-06\\
42 1.23313851646483e-05\\
43 9.42030698650241e-06\\
44 5.70123630507811e-06\\
45 4.65141152647642e-06\\
46 4.78360259825102e-06\\
47 6.50604327840325e-06\\
48 1.37435979859793e-05\\
49 0.00968102602496422\\
50 0.26697697131928\\
51 0.334537901933737\\
52 0.303791845948894\\
53 0.278179291478555\\
54 0.292025232746663\\
55 0.315172484560923\\
56 0.313620013771019\\
57 0.292630209791779\\
58 0.283028171641592\\
59 0.298141973504943\\
60 0.316637714063637\\
61 0.311128998587632\\
62 0.286783364020671\\
63 0.279965232548997\\
64 0.310260343947963\\
65 0.334509510003084\\
66 0.259051389829505\\
67 0.0140710108131047\\
68 1.69276106871839e-05\\
69 8.4323649464023e-06\\
70 6.47621244801308e-06\\
71 6.58609257308833e-06\\
72 8.64998919124006e-06\\
73 1.53255891673222e-05\\
74 4.69152086355273e-05\\
75 0.046105172775243\\
76 0.130745162553821\\
77 0.181993921121044\\
78 0.228797467450269\\
79 0.269348086901115\\
80 0.292574740049676\\
81 0.298213683901988\\
82 0.29427756042867\\
83 0.282071801496174\\
84 0.253320128760343\\
85 0.206707841552377\\
86 0.157682500038679\\
87 0.0878509565517256\\
88 2.0542992998638e-05\\
89 3.41979292892445e-06\\
90 9.96304031421011e-07\\
91 3.66625184133367e-07\\
92 1.53413264168473e-07\\
93 6.79585282907062e-08\\
94 2.93429707319927e-08\\
95 1.02312136853063e-08\\
96 6.58314272286537e-11\\
97 -5.64804231139551e-09\\
98 -8.93322282946151e-09\\
99 -1.06691268632939e-08\\
100 -1.11119661215267e-08\\
};
\node[right, inner sep=0mm, text=black]
at (axis cs:1.5,0.25,0) {NNLS};
\end{axis}

\begin{axis}[%
width=2.8in,
height=.8in,
ylabel style={yshift=-1em},
xlabel style={yshift=.5em},
scale only axis,
xmin=1,
xmax=100,
ymin=-0.01,
ymax=0.32,
at=(plot2.above north west),
anchor=below south west
]
\addplot [
color=black,
solid,
line width=2.0pt,
forget plot
]
table[row sep=crcr]{
1 0\\
2 0\\
3 0\\
4 0\\
5 0\\
6 0\\
7 0\\
8 0\\
9 0\\
10 0\\
11 0\\
12 0\\
13 0.3\\
14 0\\
15 0\\
16 0\\
17 0.2\\
18 0\\
19 0\\
20 0\\
21 0\\
22 0\\
23 0\\
24 0\\
25 0\\
26 0\\
27 0\\
28 0\\
29 0\\
30 0.3\\
31 0.3\\
32 0.3\\
33 0\\
34 0\\
35 0.3\\
36 0.3\\
37 0.3\\
38 0\\
39 0\\
40 0\\
41 0\\
42 0\\
43 0\\
44 0\\
45 0\\
46 0\\
47 0\\
48 0\\
49 0\\
50 0.3\\
51 0.3\\
52 0.3\\
53 0.3\\
54 0.3\\
55 0.3\\
56 0.3\\
57 0.3\\
58 0.3\\
59 0.3\\
60 0.3\\
61 0.3\\
62 0.3\\
63 0.3\\
64 0.3\\
65 0.3\\
66 0.3\\
67 0\\
68 0\\
69 0\\
70 0\\
71 0\\
72 0\\
73 0\\
74 0\\
75 0.05\\
76 0.1236\\
77 0.1844\\
78 0.2324\\
79 0.2676\\
80 0.29\\
81 0.2996\\
82 0.2964\\
83 0.2804\\
84 0.2516\\
85 0.21\\
86 0.1556\\
87 0.0884\\
88 0\\
89 0\\
90 0\\
91 0\\
92 0\\
93 0\\
94 0\\
95 0\\
96 0\\
97 0\\
98 0\\
99 0\\
100 0\\
};
\addplot [
color=red,
dashed,
line width=2.0pt,
forget plot
]
table[row sep=crcr]{
1 -6.12561386293141e-09\\
2 -6.66364270236864e-09\\
3 -3.62046713548723e-09\\
4 -1.39243562626484e-09\\
5 -5.15218228542785e-10\\
6 7.29827810338706e-11\\
7 7.62954446487365e-10\\
8 2.01647837588818e-09\\
9 5.26675116187286e-09\\
10 1.67029178631729e-08\\
11 8.2229156370488e-08\\
12 1.59651608727353e-06\\
13 0.299995944636277\\
14 7.36768535480115e-07\\
15 1.03700468101728e-06\\
16 7.82678879725294e-06\\
17 0.199990205054542\\
18 9.91735019616037e-07\\
19 4.47298546947789e-07\\
20 3.32209038981297e-07\\
21 3.25214960346421e-07\\
22 3.8643087710662e-07\\
23 5.45475622070738e-07\\
24 9.13512030528925e-07\\
25 1.86418554161527e-06\\
26 4.87561965240855e-06\\
27 1.78009416772787e-05\\
28 0.000102997812091545\\
29 0.000915236087559358\\
30 0.28707136558814\\
31 0.347905489462623\\
32 0.217783885838611\\
33 0.0530473867050448\\
34 0.0502230779187584\\
35 0.167131138472745\\
36 0.414975742889022\\
37 0.256198498163468\\
38 0.00336270012225539\\
39 0.000730064509108696\\
40 0.000273733301435595\\
41 0.000141887285574395\\
42 9.57399489091194e-05\\
43 8.07638137566362e-05\\
44 8.16767477581988e-05\\
45 9.92875163041986e-05\\
46 0.000148979447869595\\
47 0.000289731634150356\\
48 0.000790128276143049\\
49 0.00371699044529348\\
50 0.216796840802941\\
51 0.545250844118662\\
52 0.0644365217882123\\
53 0.186279267270328\\
54 0.68444431835943\\
55 0.155059165298259\\
56 0.114011739947842\\
57 0.324130401183207\\
58 0.57133519370097\\
59 0.191795303130506\\
60 0.126610418277379\\
61 0.306207825147452\\
62 0.592502921130692\\
63 0.141832435680519\\
64 0.0890052805901768\\
65 0.607072185816199\\
66 0.162326108048782\\
67 0.0103140524918887\\
68 0.00288080409376788\\
69 0.00138074226800169\\
70 0.000935101072978619\\
71 0.0008503181284308\\
72 0.00101650555645177\\
73 0.00162164916115143\\
74 0.00366138730735345\\
75 0.0142982065022975\\
76 0.156063904117847\\
77 0.232694143332898\\
78 0.140585537631622\\
79 0.27223135469855\\
80 0.384388293783364\\
81 0.22953870801241\\
82 0.286791787196125\\
83 0.300533650080743\\
84 0.275934336811241\\
85 0.165328936111728\\
86 0.183739696689576\\
87 0.0799176504692477\\
88 0.000984971984815399\\
89 1.81471307754817e-05\\
90 4.32562468966423e-07\\
91 2.28679206218435e-07\\
92 1.34379802025009e-07\\
93 6.65048272339596e-08\\
94 3.25465217106517e-08\\
95 1.63330122239749e-08\\
96 8.56955770249576e-09\\
97 4.74864999867258e-09\\
98 2.75949600251041e-09\\
99 1.71626063496823e-09\\
100 1.15170379006318e-09\\
};
\node[right, inner sep=0mm, text=black]
at (axis cs:1.5,0.25,0) {L1 reg};
\end{axis}

\begin{axis}[%
width=2.8in,
height=.8in,
ylabel style={yshift=-1em},
xlabel style={yshift=.5em},
scale only axis,
xmin=1,
xmax=100,
xlabel={$\text{sample number }{\it k}$},
ymin=-0.01,
ymax=0.32,
at=(plot2.below south west),
anchor=above north west
]
\addplot [
color=black,
solid,
line width=2.0pt,
forget plot
]
table[row sep=crcr]{
1 0\\
2 0\\
3 0\\
4 0\\
5 0\\
6 0\\
7 0\\
8 0\\
9 0\\
10 0\\
11 0\\
12 0\\
13 0.3\\
14 0\\
15 0\\
16 0\\
17 0.2\\
18 0\\
19 0\\
20 0\\
21 0\\
22 0\\
23 0\\
24 0\\
25 0\\
26 0\\
27 0\\
28 0\\
29 0\\
30 0.3\\
31 0.3\\
32 0.3\\
33 0\\
34 0\\
35 0.3\\
36 0.3\\
37 0.3\\
38 0\\
39 0\\
40 0\\
41 0\\
42 0\\
43 0\\
44 0\\
45 0\\
46 0\\
47 0\\
48 0\\
49 0\\
50 0.3\\
51 0.3\\
52 0.3\\
53 0.3\\
54 0.3\\
55 0.3\\
56 0.3\\
57 0.3\\
58 0.3\\
59 0.3\\
60 0.3\\
61 0.3\\
62 0.3\\
63 0.3\\
64 0.3\\
65 0.3\\
66 0.3\\
67 0\\
68 0\\
69 0\\
70 0\\
71 0\\
72 0\\
73 0\\
74 0\\
75 0.05\\
76 0.1236\\
77 0.1844\\
78 0.2324\\
79 0.2676\\
80 0.29\\
81 0.2996\\
82 0.2964\\
83 0.2804\\
84 0.2516\\
85 0.21\\
86 0.1556\\
87 0.0884\\
88 0\\
89 0\\
90 0\\
91 0\\
92 0\\
93 0\\
94 0\\
95 0\\
96 0\\
97 0\\
98 0\\
99 0\\
100 0\\
};
\addplot [
color=red,
dashed,
line width=2.0pt,
forget plot
]
table[row sep=crcr]{
1 1.35505051623852e-10\\
2 7.03746505514857e-11\\
3 2.76626499484678e-11\\
4 -1.09773301559812e-12\\
5 -1.85929605045487e-11\\
6 -2.46532239067676e-11\\
7 -1.67057478961397e-11\\
8 1.27536869953815e-11\\
9 8.50404191510279e-11\\
10 2.69570865629731e-10\\
11 8.7383755786874e-10\\
12 4.6295958733289e-09\\
13 0.29999998003236\\
14 1.51963126526766e-08\\
15 6.36519814634084e-09\\
16 8.94677332219374e-09\\
17 0.199999962490132\\
18 1.98783485227949e-08\\
19 1.90404958466672e-09\\
20 5.59652935105959e-10\\
21 2.18312534716603e-10\\
22 9.45719058620398e-11\\
23 4.71834238346958e-11\\
24 3.99432709130565e-11\\
25 6.78984091173618e-11\\
26 1.48444756487009e-10\\
27 3.45712625193784e-10\\
28 9.10269415399512e-10\\
29 4.16889178733726e-09\\
30 0.299999573664969\\
31 0.300000975835433\\
32 0.299999355800969\\
33 4.88882542226676e-08\\
34 6.80087606230728e-08\\
35 0.299999493577498\\
36 0.300000942502794\\
37 0.29999949429507\\
38 1.74422690735199e-08\\
39 4.97905394691855e-09\\
40 2.56388427155585e-09\\
41 1.71927105885317e-09\\
42 1.40945499715883e-09\\
43 1.3851235158846e-09\\
44 1.61529534192795e-09\\
45 2.21608265071893e-09\\
46 3.55998674983482e-09\\
47 6.75070016553647e-09\\
48 1.58589362198036e-08\\
49 4.88064637982344e-07\\
50 0.299998643204629\\
51 0.30000047605462\\
52 0.30000049965394\\
53 0.300000462424892\\
54 0.300000008087634\\
55 0.299998894618608\\
56 0.299999741293505\\
57 0.3000004905681\\
58 0.300000655042638\\
59 0.30000044212564\\
60 0.299999991774999\\
61 0.299997719673862\\
62 0.300000638806604\\
63 0.300000856173714\\
64 0.300000896507839\\
65 0.300000868294858\\
66 0.29999623947221\\
67 1.82293324557659e-06\\
68 9.52910940710261e-08\\
69 6.46301500650814e-08\\
70 6.59147997295051e-08\\
71 8.71886235431774e-08\\
72 1.46394023570195e-07\\
73 3.30822648497708e-07\\
74 1.47251292886841e-06\\
75 0.0499953539002308\\
76 0.123579029029116\\
77 0.184459051016396\\
78 0.232399644632406\\
79 0.267461455094846\\
80 0.290130458410481\\
81 0.29966412668933\\
82 0.296269353521253\\
83 0.280312206684662\\
84 0.251914071488882\\
85 0.20970364262649\\
86 0.155737047495716\\
87 0.0883713325140816\\
88 9.41650717234399e-07\\
89 1.03750146529524e-07\\
90 2.85800862420693e-08\\
91 1.06612564332309e-08\\
92 4.8838502131332e-09\\
93 2.63278182588422e-09\\
94 1.63438051981046e-09\\
95 1.15924320143179e-09\\
96 9.41078437399767e-10\\
97 8.79158856825768e-10\\
98 9.47693090669333e-10\\
99 1.17547038769317e-09\\
100 1.66994451600289e-09\\
};
\node[right, inner sep=0mm, text=black]
at (axis cs:1.5,0.25,0) {BOXLS};
\end{axis}
\end{tikzpicture}%

%% file: go1DSim16_reso.tex
%
%
\begin{tikzpicture}[font=\small]

\begin{axis}[%
width=2.8in,
height=1.8in,
ylabel style={yshift=-1em},
xlabel style={yshift=.5em},
scale only axis,
xmin=1,
xmax=100,
xlabel={$\text{sample number }{\it k}$},
ymin=0.5,
ymax=6.8,
ylabel={resolution (samples)},
legend style={at={(0.214732142857143,0.565476190476189)},anchor=south west,draw=black,fill=white,legend cell align=left}
]
\addplot [
color=black,
dashed,
line width=2.0pt
]
table[row sep=crcr]{
1 4.15558348704277\\
2 4.22419375269122\\
3 5.56985865432999\\
4 6.18118800626359\\
5 5.5945314915828\\
6 5.67583034119562\\
7 6.49590381359906\\
8 6.09218869017035\\
9 6.04857963842394\\
10 6.45252070668596\\
11 6.12204768407222\\
12 6.08629420618421\\
13 6.57021601358714\\
14 6.32817925169241\\
15 6.34664896313566\\
16 6.50752173537959\\
17 6.36849563100694\\
18 6.37809275098589\\
19 6.5497401853414\\
20 6.39115052857611\\
21 6.38542885198439\\
22 6.52170063752295\\
23 6.40694406472934\\
24 6.39003489593787\\
25 6.53385634451172\\
26 6.41393492050941\\
27 6.42818338187512\\
28 6.52073121180897\\
29 6.45091427013114\\
30 6.43596662084497\\
31 6.52442170782968\\
32 6.45774613225316\\
33 6.43511421048578\\
34 6.51676484383332\\
35 6.46145551070544\\
36 6.4474617775089\\
37 6.51786659027486\\
38 6.46419618971362\\
39 6.44667658732063\\
40 6.51240016432069\\
41 6.46305507695518\\
42 6.44871009083281\\
43 6.51014866139785\\
44 6.47409049839867\\
45 6.46234070920827\\
46 6.50571723677998\\
47 6.4760603665773\\
48 6.46222264836803\\
49 6.50637143881149\\
50 6.47898399280325\\
51 6.44983404160502\\
52 6.50260405288154\\
53 6.48077864170797\\
54 6.44986343159721\\
55 6.50228602255061\\
56 6.48316749127169\\
57 6.44745528785958\\
58 6.4967110534508\\
59 6.47683797373539\\
60 6.44659371364333\\
61 6.49554609953001\\
62 6.4793027518063\\
63 6.44298389333942\\
64 6.49032633564147\\
65 6.48144543683015\\
66 6.44235065463393\\
67 6.48848686930983\\
68 6.48408866197097\\
69 6.44102338668959\\
70 6.4836880528829\\
71 6.46867419747705\\
72 6.40792745326436\\
73 6.47367055951869\\
74 6.47074601247406\\
75 6.39758363493324\\
76 6.46740495901467\\
77 6.47089245321111\\
78 6.37788710707059\\
79 6.46123230281193\\
80 6.42679448897075\\
81 6.36392527680104\\
82 6.42872957636927\\
83 6.4314262395629\\
84 6.36012260251315\\
85 6.41552051387727\\
86 6.43000341028456\\
87 6.08284638345473\\
88 6.38477242295781\\
89 6.4186135994487\\
90 6.04189624623449\\
91 6.28912966725999\\
92 6.13965593875764\\
93 6.03378415643259\\
94 6.16545503882582\\
95 5.96805738258515\\
96 5.21221777198948\\
97 6.11670715752072\\
98 5.4760371082934\\
99 4.66016912862933\\
100 3.54086101813854\\
};
\addlegendentry{EC case};

\addplot [
color=blue,
solid,
line width=2.0pt
]
table[row sep=crcr]{
1 1.00000000014801\\
2 1.0000000001631\\
3 1.00000000017659\\
4 1.00000000024949\\
5 1.00000000018395\\
6 1.00000000017675\\
7 1.00000000018563\\
8 1.00000000021383\\
9 1.00000000018419\\
10 1.00000000021237\\
11 1.00000000018025\\
12 1.000000000286\\
13 1.00000000029351\\
14 1.00000000017262\\
15 1.00000000027662\\
16 1.00000000033347\\
17 1.00000000022793\\
18 1.00000000024105\\
19 1.00000000018543\\
20 1.00000000022035\\
21 1.00000000025463\\
22 1.00000000039245\\
23 1.00000000024679\\
24 1.00000000021812\\
25 1.0000000002746\\
26 1.00000000022263\\
27 1.00000000024685\\
28 1.00000000022053\\
29 1.00000000020538\\
30 1.20362020861563\\
31 1.73449453273171\\
32 2.31988958820425\\
33 2.24568734213135\\
34 2.29912731760583\\
35 2.45767888844893\\
36 2.4039191280779\\
37 2.21163640996834\\
38 1.41875069602875\\
39 1.00000000019277\\
40 1.00000000022081\\
41 1.00000000030077\\
42 1.00000000022015\\
43 1.00000000024382\\
44 1.00000000037101\\
45 1.00000000026468\\
46 1.0000000002215\\
47 1.00000000027717\\
48 1.03065171285491\\
49 1.89124944251614\\
50 2.69694911809852\\
51 2.89984109908763\\
52 3.16106311580145\\
53 3.56043692402571\\
54 3.58327613495115\\
55 3.81282009335379\\
56 3.86407188722862\\
57 3.94734612354095\\
58 3.98886977093781\\
59 3.94150465256492\\
60 4.00408920589558\\
61 3.978446880005\\
62 4.00293500000688\\
63 3.83653794081916\\
64 3.76784777573108\\
65 3.77728384650612\\
66 3.54279928878087\\
67 3.37582619919123\\
68 3.07150717830697\\
69 2.80588391047112\\
70 2.61516795629206\\
71 2.50965356848394\\
72 2.73348101727667\\
73 2.85755022435228\\
74 3.17126514306303\\
75 3.26157814342081\\
76 3.51122501958529\\
77 3.55901440154213\\
78 3.5635318733764\\
79 3.30872956789062\\
80 3.49703186020318\\
81 3.28700844352432\\
82 3.27622714205737\\
83 3.08430078507369\\
84 2.59329856985287\\
85 2.5597609003888\\
86 2.13415657935509\\
87 1.3003882550226\\
88 1.00000000020099\\
89 1.00000000033158\\
90 1.00000000027907\\
91 1.00000000024874\\
92 1.00000000030861\\
93 1.00000000021606\\
94 1.00000000018795\\
95 1.00000000024299\\
96 1.00000000015293\\
97 1.0000000002794\\
98 1.00000000027791\\
99 1.00000000015116\\
100 1.00000000020379\\
};
\addlegendentry{NN case};

\addplot [
color=red,
dash pattern=on 1pt off 3pt on 3pt off 3pt,
line width=2.0pt
]
table[row sep=crcr]{
1 1.00000000119005\\
2 1.00000000220286\\
3 1.00000002586107\\
4 1.00000000005126\\
5 1.00000000649056\\
6 1.00000001414892\\
7 0.999999999985094\\
8 1.00000000000039\\
9 1.00000000000238\\
10 1.0000000104712\\
11 1.0000000032313\\
12 1.00000000268082\\
13 1.00000000967048\\
14 1.00000000000908\\
15 1.00000000027261\\
16 1.00000000328604\\
17 1.00000000601322\\
18 1.00000001009104\\
19 1.00000000112097\\
20 1.00000000348694\\
21 1.00000000074663\\
22 1.00000002304512\\
23 1.00000000138398\\
24 1.00000000000037\\
25 1.00000000444781\\
26 1.0000000023796\\
27 0.999999999998353\\
28 1.00000000389895\\
29 1.00000000334875\\
30 1.00000000323175\\
31 1.00000000243743\\
32 1.0000000043029\\
33 1.00000000036895\\
34 1.00000001137163\\
35 1.00000001165893\\
36 1.00000000008867\\
37 1.00000000605001\\
38 1.00000000771052\\
39 1.00000000411763\\
40 1.00000000001532\\
41 1.00000000001258\\
42 1.00000000004876\\
43 1.00000000039312\\
44 1.00000000286815\\
45 1.00000000533178\\
46 1.00000000050211\\
47 1.00000000113492\\
48 1.00000001956004\\
49 1.00000000924082\\
50 1.0000000039291\\
51 1.00000000043719\\
52 1.00000000004297\\
53 1.0000000023422\\
54 1.00000000704067\\
55 1.00000000001096\\
56 1.00000002313984\\
57 1.0000000000457\\
58 1.00000001833394\\
59 1.00000000000285\\
60 1.00000000917603\\
61 1.00000000366263\\
62 1.00000000000671\\
63 1.00000000837905\\
64 1.00000000378045\\
65 1.00000000750136\\
66 1.00000000199055\\
67 1.00000000378256\\
68 1.00000000000915\\
69 0.99999999997449\\
70 1.00000000300914\\
71 0.999999998975307\\
72 1.00000000166604\\
73 1.00000000617333\\
74 1.00000000510113\\
75 1.25462707511465\\
76 1.87095576587457\\
77 2.16266767610119\\
78 2.16617623018025\\
79 2.18491207756498\\
80 2.16271565994174\\
81 1.85387002700661\\
82 1.89782411862613\\
83 2.18126579105672\\
84 2.10195472806462\\
85 2.11743321597124\\
86 1.69578348452936\\
87 1.18600206940266\\
88 1.00000000493759\\
89 1.00000000318626\\
90 1.00000000566871\\
91 1.00000000635374\\
92 1.00000001488251\\
93 1.00000000001701\\
94 1.00000000054499\\
95 1.00000000340359\\
96 1.00000000005166\\
97 1.00000000030916\\
98 1.00000000250882\\
99 1.00000000336428\\
100 1.00000000009803\\
};
\addlegendentry{BOX case};

\end{axis}
\end{tikzpicture}%

%% file: go1DSim16_lowres.tex
%
%
\begin{tikzpicture}[font=\small]

\begin{axis}[%
width=2.8in,
height=.8in,
ylabel style={yshift=-1em},
scale only axis,
xmin=1,
xmax=100,
ymin=-0.01,
ymax=0.32,
ylabel={sample value},
name=plot2
]
\addplot [
color=black,
solid,
line width=2.0pt,
forget plot
]
table[row sep=crcr]{
1 0.00495727428642567\\
2 0.00495837708695035\\
3 0.00495418915625123\\
4 0.00505330180749297\\
5 0.00496361031764536\\
6 0.00496737964476779\\
7 0.00496018568446743\\
8 0.00496701705742453\\
9 0.004957131151059\\
10 0.00492271552866441\\
11 0.00495892297863065\\
12 0.00498157525726128\\
13 0.304960193288935\\
14 0.00495446682907641\\
15 0.00481823633890599\\
16 0.00496048803324811\\
17 0.204957201520301\\
18 0.0049610898240644\\
19 0.00495538607356139\\
20 0.00495653433495136\\
21 0.00496152853884269\\
22 0.00496015162934782\\
23 0.00495974039768043\\
24 0.00495368323026923\\
25 0.00495768185646739\\
26 0.00496243186989886\\
27 0.00496385651968012\\
28 0.00495602611408685\\
29 0.00513496378698619\\
30 0.301498024520697\\
31 0.305182956801218\\
32 0.244440482900245\\
33 0.122899615860661\\
34 0.0947034032615193\\
35 0.224434345647751\\
36 0.304195675872052\\
37 0.194838412749959\\
38 0.0905723861887964\\
39 0.00783235284870898\\
40 0.00534663694998017\\
41 0.00502160451651434\\
42 0.00499890268110903\\
43 0.00499343224873883\\
44 0.00498930735557224\\
45 0.00500645517058729\\
46 0.00511437762179412\\
47 0.00590243456827011\\
48 0.0112217199348379\\
49 0.0557681971495185\\
50 0.196818176882516\\
51 0.308550871915941\\
52 0.302984324924182\\
53 0.306424341455568\\
54 0.304776042716185\\
55 0.306618255666308\\
56 0.305635070977587\\
57 0.304867262319021\\
58 0.303229512730468\\
59 0.30377191757907\\
60 0.307232380366258\\
61 0.306860086855522\\
62 0.304485211092469\\
63 0.306968117338329\\
64 0.304978538872092\\
65 0.2612134463252\\
66 0.186622225432075\\
67 0.0880726272052925\\
68 0.0336054747167509\\
69 0.0120106053072959\\
70 0.0097459842918397\\
71 0.0097430777086629\\
72 0.00955817103385925\\
73 0.0123510830162559\\
74 0.0337524754277183\\
75 0.0707591193981898\\
76 0.124916722514172\\
77 0.180467323905759\\
78 0.223156952932186\\
79 0.257960406419443\\
80 0.288869811635095\\
81 0.295515131880165\\
82 0.28768028387276\\
83 0.274639664093115\\
84 0.258845188966006\\
85 0.199553963033395\\
86 0.168860243018571\\
87 0.113090729867508\\
88 0.00697824830422178\\
89 0.00500926020094994\\
90 0.00496740207745461\\
91 0.00495959041472815\\
92 0.00493110856041312\\
93 0.00496010132337688\\
94 0.00495532651439134\\
95 0.00495846945705125\\
96 0.00495746926753782\\
97 0.00494829801027663\\
98 0.00496268165352376\\
99 0.00496334182162173\\
100 0.00496032844534966\\
};
\addplot [
color=black,
solid,
line width=2.0pt,
forget plot
]
table[row sep=crcr]{
1 -0.00495502525518532\\
2 -0.00495646399031102\\
3 -0.00496309869731704\\
4 -0.00486943837313447\\
5 -0.00496117304828658\\
6 -0.00495179884455865\\
7 -0.00495631669036811\\
8 -0.00495380212578311\\
9 -0.00496072244754941\\
10 -0.00500719393676263\\
11 -0.00496123979337426\\
12 -0.00494275679557177\\
13 0.295040558995879\\
14 -0.00495326878444757\\
15 -0.00510212020890322\\
16 -0.00495899139127687\\
17 0.195044537884087\\
18 -0.00495654028418357\\
19 -0.00495953879908484\\
20 -0.00495842187638118\\
21 -0.00495342769772833\\
22 -0.00495941165172553\\
23 -0.00495665806010948\\
24 -0.00495645222781604\\
25 -0.00496679503157793\\
26 -0.00495601103648369\\
27 -0.00495566757604138\\
28 -0.00495924666688552\\
29 -0.00479590009672393\\
30 0.291498024638713\\
31 0.295182956888311\\
32 0.23444048381316\\
33 0.11289961620605\\
34 0.0847034050184448\\
35 0.214434349225485\\
36 0.294195676566687\\
37 0.184838413000193\\
38 0.0805723862040395\\
39 -0.00214506236125089\\
40 -0.00457216936683835\\
41 -0.00490104139271352\\
42 -0.0049156947152369\\
43 -0.00491492634500901\\
44 -0.00492894814465217\\
45 -0.00491607041340103\\
46 -0.00480547331244452\\
47 -0.00403084293157008\\
48 0.00122172000919818\\
49 0.0457681971593047\\
50 0.18681817794004\\
51 0.298550871766565\\
52 0.292984325653379\\
53 0.296424342464888\\
54 0.294776048549465\\
55 0.296618256110378\\
56 0.295635072492587\\
57 0.294867263392916\\
58 0.29322951392146\\
59 0.293771918695938\\
60 0.297232386042197\\
61 0.296860088939866\\
62 0.294485211487427\\
63 0.296968117427696\\
64 0.294978539332078\\
65 0.251213446756083\\
66 0.176622225553729\\
67 0.078072627288293\\
68 0.0236054724118731\\
69 0.00201060593826696\\
70 -0.000254014102665678\\
71 -0.000256921621257789\\
72 -0.000441828961356805\\
73 0.00235108297238185\\
74 0.0237524756265657\\
75 0.0607591697772705\\
76 0.114916722624912\\
77 0.170467323959201\\
78 0.213156953030079\\
79 0.247960406508355\\
80 0.278869811008008\\
81 0.285515132137789\\
82 0.277680284711096\\
83 0.264639664667357\\
84 0.248845189008762\\
85 0.189553977037576\\
86 0.158860256384287\\
87 0.103090729918222\\
88 -0.00298213061523711\\
89 -0.00490628173645291\\
90 -0.00494771753255918\\
91 -0.00495412519376259\\
92 -0.00498360674828291\\
93 -0.00495925321774848\\
94 -0.00495684603902191\\
95 -0.00495622757625824\\
96 -0.00495776355091948\\
97 -0.00496672209101234\\
98 -0.00495480016797956\\
99 -0.00495082409827319\\
100 -0.00495636528830801\\
};
\node[right, inner sep=0mm, text=black]
at (axis cs:3,0.25,0) {NN};
\end{axis}

\begin{axis}[%
width=2.8in,
height=.8in,
ylabel style={yshift=-2em},
scale only axis,
xmin=1,
xmax=100,
ymin=-0.01,
ymax=0.32,
at=(plot2.above north west),
anchor=below south west
]
\addplot [
color=black,
solid,
line width=2.0pt,
forget plot
]
table[row sep=crcr]{
1 0.000829968727870778\\
2 0.000752559211317561\\
3 0.00017897976296366\\
4 0.000501941872499856\\
5 0.000196593520811916\\
6 0.000241939811616909\\
7 0.00161437447426849\\
8 0.0076203491835507\\
9 0.0125270817129977\\
10 0.0196586894044444\\
11 0.0423208225569063\\
12 0.0492156298860138\\
13 0.0528737030044866\\
14 0.0610694137672798\\
15 0.06048364143588\\
16 0.0560280818688472\\
17 0.0438276185519122\\
18 0.036294035776403\\
19 0.0287797271650767\\
20 0.0167196453434073\\
21 0.0109207228017918\\
22 0.00652399947461874\\
23 0.00197395033890135\\
24 0.00231181938787017\\
25 0.00711556333245511\\
26 0.0265787489405369\\
27 0.0433904831468865\\
28 0.0661903743206533\\
29 0.109525566250282\\
30 0.131100993518653\\
31 0.150317895597155\\
32 0.170893274031374\\
33 0.176947011599161\\
34 0.176943503665385\\
35 0.169180589468516\\
36 0.156066774239134\\
37 0.137526477977793\\
38 0.100939787008979\\
39 0.0758261244314085\\
40 0.0524655497899822\\
41 0.023379906742381\\
42 0.0108574050593724\\
43 0.00437781918742708\\
44 0.0049284875834843\\
45 0.0104974584581221\\
46 0.0228489299488785\\
47 0.0531190210719146\\
48 0.0825177927603473\\
49 0.118478980613203\\
50 0.178851688566919\\
51 0.213946085303508\\
52 0.246058255344386\\
53 0.278003544982191\\
54 0.289449504844784\\
55 0.294471952513905\\
56 0.296856166443362\\
57 0.29646794553874\\
58 0.295453233960927\\
59 0.296114407681786\\
60 0.297170300891881\\
61 0.295991381323127\\
62 0.288563701035603\\
63 0.276974197009384\\
64 0.256653509504098\\
65 0.210781245821811\\
66 0.175262086971866\\
67 0.136516252068986\\
68 0.0785750341698468\\
69 0.0500941154451731\\
70 0.0301359604805675\\
71 0.0167786511966046\\
72 0.0185260067212099\\
73 0.0308015907411701\\
74 0.0624681550874318\\
75 0.0861211866566478\\
76 0.11816464866299\\
77 0.172272792690824\\
78 0.199742639090182\\
79 0.224903834740145\\
80 0.253874679232666\\
81 0.261569596237933\\
82 0.260204955886691\\
83 0.237286049219305\\
84 0.216508013159213\\
85 0.192393323928361\\
86 0.128039874572046\\
87 0.111014123235407\\
88 0.0790640139310064\\
89 0.0311587519609242\\
90 0.0216087005836485\\
91 0.0138566180300346\\
92 0.00230231221984045\\
93 0.00144953750161037\\
94 0.00260239694882003\\
95 0.000620671790926508\\
96 0.00191181871574742\\
97 0.000791099517776672\\
98 0.000969783168697447\\
99 0.000192509592386259\\
100 0.00100823968318764\\
};
\node[right, inner sep=0mm, text=black]
at (axis cs:3,0.25,0) {EC};
\end{axis}

\begin{axis}[%
width=2.8in,
height=.8in,
ylabel style={yshift=-2em},
xlabel style={yshift=.5em},
scale only axis,
xmin=1,
xmax=100,
xlabel={$\text{sample number }{\it k}$},
ymin=-0.01,
ymax=0.32,
at=(plot2.below south west),
anchor=above north west
]
\addplot [
color=black,
solid,
line width=2.0pt,
forget plot
]
table[row sep=crcr]{
1 0.00499519320160857\\
2 0.0049862702685175\\
3 0.0049850237332123\\
4 0.00499285512921688\\
5 0.00498696272089205\\
6 0.00499547369642528\\
7 0.00500344868342495\\
8 0.00498506580152025\\
9 0.00498100727179462\\
10 0.00498834432966078\\
11 0.00497695280799348\\
12 0.00499231841416758\\
13 0.304947298732799\\
14 0.00498990907354369\\
15 0.00499643717546405\\
16 0.0050007183170635\\
17 0.204933093098504\\
18 0.00499862030170561\\
19 0.00498316221279538\\
20 0.00498727735316606\\
21 0.00497477635423138\\
22 0.00499397209642893\\
23 0.00497822741110099\\
24 0.00499289514016255\\
25 0.00499388808447065\\
26 0.00499283968184727\\
27 0.00491866390245832\\
28 0.00496802248846961\\
29 0.00497296199091579\\
30 0.304968254510271\\
31 0.304985108554291\\
32 0.304980789939833\\
33 0.00501960215865438\\
34 0.00501632077111935\\
35 0.30497565822958\\
36 0.304991474546455\\
37 0.304980788833433\\
38 0.00498992727378322\\
39 0.0049812128450526\\
40 0.00499502774553662\\
41 0.00508847068962837\\
42 0.00498770138045757\\
43 0.00501641230138716\\
44 0.00498727491140016\\
45 0.00499153001541686\\
46 0.00498274475077665\\
47 0.00499457904976452\\
48 0.0049879858445081\\
49 0.00498568315072134\\
50 0.30498288285051\\
51 0.304992590078257\\
52 0.304970729217132\\
53 0.304990666471753\\
54 0.304987880044422\\
55 0.304967639837713\\
56 0.304974445631643\\
57 0.30498106870391\\
58 0.304984380941306\\
59 0.304994664843207\\
60 0.304993198387393\\
61 0.304994915750001\\
62 0.305003231829119\\
63 0.304990973650689\\
64 0.304985330500614\\
65 0.304988560729271\\
66 0.304989519677861\\
67 0.00498973572410932\\
68 0.00499388344138652\\
69 0.00500458903238155\\
70 0.00499947088604813\\
71 0.00495520536827598\\
72 0.00500519327081861\\
73 0.00516767768530713\\
74 0.00738030589644723\\
75 0.0660673487318206\\
76 0.131879874589629\\
77 0.181421972657489\\
78 0.225883729790468\\
79 0.274797255025703\\
80 0.286207717933451\\
81 0.294287742794282\\
82 0.290438922503668\\
83 0.285094375805372\\
84 0.246764942801974\\
85 0.209614254024672\\
86 0.17378051765543\\
87 0.0991452434944128\\
88 0.00663137511614309\\
89 0.00505533190801088\\
90 0.00499709826613071\\
91 0.00499729065072074\\
92 0.00499053325495424\\
93 0.00498821685843609\\
94 0.00498632949432931\\
95 0.00499147584241655\\
96 0.00503765749081708\\
97 0.00498905538901795\\
98 0.00498314721278348\\
99 0.00499603845864272\\
100 0.00500139971671842\\
};
\addplot [
color=black,
solid,
line width=2.0pt,
forget plot
]
table[row sep=crcr]{
1 -0.00498271671207817\\
2 -0.00499208884468771\\
3 -0.00499324947057289\\
4 -0.00498553229743948\\
5 -0.00499106010786932\\
6 -0.00498292715525395\\
7 -0.00497482591100606\\
8 -0.00499294898480684\\
9 -0.00499723822500187\\
10 -0.00498971517015434\\
11 -0.0050018879937852\\
12 -0.00498570190558212\\
13 0.29496909120661\\
14 -0.0049885120344868\\
15 -0.00498171394170299\\
16 -0.00497731677799607\\
17 0.194953719165795\\
18 -0.00497933469029022\\
19 -0.00499542684099197\\
20 -0.00499056207146964\\
21 -0.00500394190655129\\
22 -0.00498426572403332\\
23 -0.00500122816865201\\
24 -0.00498533047175442\\
25 -0.00498496765780487\\
26 -0.00498557931450705\\
27 -0.00505967633552018\\
28 -0.00501007956060562\\
29 -0.00500508769147245\\
30 0.294990619512276\\
31 0.295007568884216\\
32 0.295003812620138\\
33 -0.00495880783768143\\
34 -0.00496170839772958\\
35 0.294996772236772\\
36 0.2950132444536\\
37 0.295002681221234\\
38 -0.00498706223936551\\
39 -0.00499643695945906\\
40 -0.00498394449709849\\
41 -0.00489103512589395\\
42 -0.00499001980631419\\
43 -0.00496194771534419\\
44 -0.0049906089845706\\
45 -0.00498682569559605\\
46 -0.00499503633912468\\
47 -0.00498293703682862\\
48 -0.00498960632648959\\
49 -0.00499256382926205\\
50 0.295005588222239\\
51 0.295014140796638\\
52 0.294992635765936\\
53 0.295013161184272\\
54 0.295009968308165\\
55 0.294988954889163\\
56 0.294995463566778\\
57 0.295003239348461\\
58 0.295007067812094\\
59 0.295015900447169\\
60 0.295016215421413\\
61 0.295016396739385\\
62 0.295024338129508\\
63 0.295013372156802\\
64 0.295006103027248\\
65 0.29500931077996\\
66 0.295012152300217\\
67 -0.0049884893397234\\
68 -0.00498560387239877\\
69 -0.0049733613713272\\
70 -0.00497825045670197\\
71 -0.00502366454139747\\
72 -0.00497247666311296\\
73 -0.00481062833077317\\
74 -0.00260570196144272\\
75 0.0560671847985077\\
76 0.121879225636036\\
77 0.17142492181614\\
78 0.21588334210216\\
79 0.264805577142649\\
80 0.276207658978969\\
81 0.284287673655427\\
82 0.280438675068638\\
83 0.275094397441114\\
84 0.236765004761622\\
85 0.199617886667795\\
86 0.163780802457012\\
87 0.089124586844747\\
88 -0.00335255654140099\\
89 -0.0049226848198316\\
90 -0.0049812531808584\\
91 -0.00498058185991113\\
92 -0.00498686640645474\\
93 -0.0049910923163452\\
94 -0.00499212401328464\\
95 -0.00498625112504669\\
96 -0.00493997119019696\\
97 -0.00498875166465496\\
98 -0.00499564929330631\\
99 -0.0049814964595214\\
100 -0.00497686511306483\\
};
\node[right, inner sep=0mm, text=black]
at (axis cs:3,0.25,0) {BOX};
\end{axis}
\end{tikzpicture}%